%% file: main.tex
\documentclass[12pt,a4paper,leqno]{article}

\usepackage{latexsym}
\usepackage{amssymb,exscale}
\usepackage[centertags]{amsmath}
\usepackage{amsthm}
\usepackage[all]{xy}
\usepackage{graphicx}

\input{preamble}

\begin{document}

\input{intro}

\input{sect1}
\input{sect2}
\input{sect3}

\input{references}
\bigskip
\filbreak
{
\small
\begin{sc}
\parindent=0pt\baselineskip=12pt
\parbox{2.3in}{
Mikhail Sodin\\
School of Mathematics\\
Tel Aviv University\\
Tel Aviv 69978, Israel
\smallskip
\emailwww{sodin@tau.ac.il}
{}
}
\hfill
\parbox{2.3in}{
Boris Tsirelson\\
School of Mathematics\\
Tel Aviv University\\
Tel Aviv 69978, Israel
\smallskip
\emailwww{tsirel@tau.ac.il}
{www.tau.ac.il/\textasciitilde tsirel/}
}
\end{sc}
}
\filbreak

\end{document}

%% file: preamble.tex
\numberwithin{equation}{section}

\swapnumbers
\newtheorem{theorem}[equation]{Theorem}
\theoremstyle{definition}

\newtheorem{definition}[equation]{Definition}

\newtheorem{example}[equation]{Example}

\theoremstyle{plain}

\renewcommand{\phi}{\varphi}

\renewcommand{\(}{\bigl(}
\renewcommand{\)}{\bigr)\vphantom{)}}

\newcommand{\ip}[2]{\langle#1,#2\rangle}

\newcommand{\const}{\operatorname{const}}

\newcommand{\Var}{\operatorname{Var}}
\newcommand{\dist}{\operatorname{dist}}
\renewcommand{\Re}{\operatorname{Re}}
\renewcommand{\Im}{\operatorname{Im}}

\newcommand{\Om}{\Omega}
\newcommand{\om}{\omega}

\newcommand{\De}{\Delta}
\newcommand{\cN}{\mathcal N}
\newcommand{\cL}{\mathcal L}
\newcommand{\M}{\mathcal M}

\newcommand{\Ex}{\mathbb E}
\renewcommand{\Pr}[1]{\,\mathbb P\,\(\,#1\,\)\,}
\newcommand{\R}{\mathbb R}
\newcommand{\C}{\mathbb C}
\newcommand{\Z}{\mathbb Z}
\newcommand{\N}{\mathbb N}
\renewcommand{\S}{\mathbb S}

\newcommand{\D}{\mathbb D}

\newcommand{\ti}{\tilde}

\newcommand{\Wi}[1]{:#1:}

\def\emailwww#1#2{\par\quad {\tt #1}\par\quad {\tt #2}\medskip}

%% file: intro.tex
\title{Random complex zeroes, I. \\
 Asymptotic normality}
\author{Mikhail Sodin\thanks{Supported by the Israel
Science Foundation of the Israel Academy
of Sciences and Humanities} 
\ and Boris Tsirelson}

\date{}

\maketitle

\begin{abstract}
We consider three models (elliptic, flat and hyperbolic) of Gaussian
random analytic functions distinguished by invariance of their zeroes
distribution. Asymptotic normality is proven for smooth functionals
(linear statistics) of the set of zeroes.
\end{abstract}

\section*{Introduction and the main result}

Zeroes of random polynomials and other analytic functions were
studied by mathematicians and physicists under various assumptions
on random coefficients. One class of models introduced not long
ago by Bogomolny, Bohigas and Leboeuf \cite{BBL92, BBL96},
Kostlan \cite{Ko93}, and Shub and Smale \cite{ShSm} has a
remarkably unique unitary invariance:
\begin{quote}
``\dots indeed it has no true freedom at all. It is (statistically)
unique in the same sense as `the Poisson process', or `the thermal
(black body) electromagnetic field' are unique\dots''\\
Hannay \cite[p.~L755]{Ha98}
\end{quote}
Following Hannay \cite{Ha96}, we use the term `chaotic analytic
zero points' (CAZP, for short). We consider here three CAZP
models: the elliptic CAZP, the flat CAZP, and the hyperbolic CAZP
called by Leboeuf \cite[p.~654]{Le99} $ SU (2) $, $ W_1 $, and $
SU (1,1) $ in accordance with the symmetry group of the model.
These models may be described analytically or geometrically. The
analytical description is short and elementary: CAZP is the
(random) set of zeroes of such a (random) analytic function $ \psi
$,
\begin{align}
\psi (z) &= \sum_{k=0}^L \zeta_k \sqrt{
 \frac{L(L-1) \dots (L-k+1)}{k!} } \, z^k
 && \text{(elliptic, $ L = 1, 2, \dots $)},
 \label{elliptic} \displaybreak[0] \\
\psi (z) &= \sum_{k=0}^\infty \zeta_k \sqrt{
 \frac{L^{k}}{k!}} \, z^k
 && \text{(flat, $ L>0 $)}, \label{flat} \displaybreak[0] \\
\psi (z) &= \sum_{k=0}^\infty \zeta_k \sqrt{
 \frac{L(L+1) \dots(L+k-1)}{k!}}\, z^k
 && \text{(hyperbolic, $ L>0$)},
\label{hyperbolic}
\end{align}
where $ \zeta_0, \zeta_1, \dots $ are independent standard complex
Gaussian (random) variables; that is, the distribution $
\cN_\C(0,1) $ of each $ \zeta_k $ has the density $ \pi^{-1}
\exp(-|z|^2) $ with respect to the Lebesgue measure $ m $ on $
\C $, $ dm(z) = (d\text{Re} z) (d\text{Im} z) $.
For the elliptic CAZP one could assume equivalently that $
(\zeta_0,\dots,\zeta_L) $ is uniformly distributed on the sphere $
|\zeta_0|^2 + \dots + |\zeta_L|^2 = \text{const} $ (which changes
the function $ \psi $ without changing its zeroes). The analytic
function (\ref{elliptic}) is a polynomial of degree $ L $, the
function (\ref{flat}) with probability one is an entire function
(indeed, $ \limsup |\zeta_k|^{1/k} = 1 $ a.s.),
and the function (\ref{hyperbolic}) with probability one is
analytic in the unit disk.

Why just these coefficients? Because of symmetry of the models
revealed by a geometric description, given in Sect.~1; readers not
bothered by this question may skip that section.

We introduce unified notations: $ \M $ for the domain of the random
function, and $ G $ for the symmetry group.
\[
\renewcommand{\arraystretch}{1.9}
\begin{array}{|c|c|c|c|} \hline
 & \M \text{ \small (domain)} &
  \begin{smallmatrix} \text{\small invariant} \\
                      \text{\small metric} \end{smallmatrix}
 & G \text{ \small (symmetries)} \\ \hline
\text{\small Elliptic} & 
  \C \cup \{\infty\} \text{\small, that is, } \S^2
 & \dfrac{ |dz| }{ 1+|z|^2 } &
 z \mapsto \dfrac{ az+b }{ -\bar b z+\bar a }, \, |a|^2+|b|^2=1 \\
\text{\small Flat} & \C & |dz| & z \mapsto az+b, \, |a|=1 \\
\text{\small Hyperbolic} & \D = \{ z \in \C : |z| < 1 \} & \dfrac{ |dz| }{ 1-|z|^2 }
  & z \mapsto \dfrac{az+b}{\bar b z+\bar a}, \, |a|^2-|b|^2=1 \\ \hline
\end{array}
\]

A symmetry $ g \in G $ transforms the random function $ \psi $ into
another random function $ z \mapsto \psi(g(z)) $ whose distribution
depends on $ g $. However, for every $ g \in G $ there exists a
complex-valued function $ u_g $ on $ \M $, $ |u_g(z)| = 1 $,
such that two random functions
\begin{equation}\label{eq0.35}
z \mapsto \frac{ \psi (g(z)) }{ \| \psi (g(z)) \| } \quad \text{and}
\quad z \mapsto u_g(z) \frac{ \psi (z) }{ \| \psi(z) \| }
\end{equation}
are identically distributed. Here
\[
\| \psi (z) \| = \left( \Ex |\psi(z)|^2 \right)^{1/2} =
\begin{cases}
 (1+|z|^2)^{L/2} &\text{(elliptic)},\\
 \exp (L|z|^2/2) &\text{(flat)}, \\
 (1-|z|^2)^{-L/2} &\text{(hyperbolic)}.
\end{cases}
\]
For example, a shift $ g(z) = z+b $ (for the flat
case); here $ u_g(z) = \exp ( iL \Im (z\bar b) ) $. The
symmetry ensures that \emph{the distribution of the (random) set of
zeroes of $ \psi $ is invariant under $ G $}. For details, see
Sect.~1.

In each of the three cases we have a parameter $ L $ which is
assumed to be large. Increasing $ L $, we increase accordingly the
mean number of random zeroes per unit area of $ \M $. 

In the flat case, there is
another interpretation: introducing $ L $, we just
make a homothety of the plane with coefficient $ \sqrt L $.
This makes some computations simpler. 

In the other two models, changing the value of the parameter $ L $,
we can change the properties of the process. We mention a recent discovery 
of Peres and Vir\'ag \cite{PV} that in the special case $ L=1 $ of the 
hyperbolic model, the point process is a {\em determinantal} one 
\cite{Soshnikov}. In this case, 
Peres and Vir\'ag found a simple expression for $n$-point correlation
functions and explicitly computed the distribution of random zeroes
in discs in $ \mathbb D $. 
        
Below, the asymptotic normality is stated for smooth functionals 
of the random set of zeroes, namely, for the random variables (often
called the \emph{linear statistics of zeroes})
\begin{equation}
Z_L (h) = \sum_{ z : \psi_L (z)=0 } h(z) = \int_{\M} h \,
dn_{\psi_L}
\end{equation}
where $h$ is a smooth test function, and $ n_{\psi} $ is the
counting measure on the set of zeroes. In the elliptic case, $
n_{\psi_L}(\S^2)=L $; in the other two cases $ n_\psi $ is
infinite but locally finite. The choice of the test-functions
depends on the model: in the elliptic case $ h $ is a real-valued
$ C^2 $-function on the Riemann sphere, in the flat case $ h $ is
a real-valued $ C^2 $-function with a compact support in the
plane, and in the hyperbolic case $ h $ is a $ C^2 $-function with
a compact support in the unit disk. Multiple zeroes may be
ignored, as well as a possible zero at $ \infty $ in the elliptic
case, since it happens with zero probability.

\medskip\noindent\textsc{Expectation of the linear statistics.}
For the elliptic model, invariance itself gives us the expectation of
the random variable $ Z_L(h) $; in the other two cases it defines the
expectation up to a numerical coefficient which can be easily found using
the Edelman-Kostlan formula $ \Ex n_\psi =
(2\pi)^{-1} \Delta \log \| \psi \| \, dm $ (see \cite{EK}, \cite{So}):
\[
\Ex Z_L(h)  = L \frac1{\pi} \int_{\M} h \, dm^*
\]
where
\[
dm^*(z) = \begin{cases}
dm(z) &\text{(flat)}, \\
(1+|z|^2)^{-2} \, dm(z) &\text{(elliptic)}, \\
(1-|z|^2)^{-2} \, dm(z) &\text{(hyperbolic)}
\end{cases}
\]
is the invariant measure on $ \M $.

\medskip\noindent\textsc{Fluctuations of the linear statistics.}
These are computed in Section~2. We get
\begin{equation}\label{0.8L}
\Var Z_L (h) = \frac{ \kappa }{ L } \| \Delta^* h
\|^2_{L^2(m^*)} + o(L^{-1}) \, , \qquad  L\to\infty \, ,
\end{equation}
where $ \kappa $ is a numerical constant (the same for each of the
three cases), and  $ \Delta^* $ is the invariant Laplacian on $
\M $:
\[
\Delta^* h(z) = \begin{cases}
\Delta h(z) & \text{(flat)}, \\
(1+|z|^2)^2 \Delta h (z) & \text{(elliptic)}, \\
(1-|z|^2)^2 \Delta h(z)  &\text{(hyperbolic)}.
\end{cases}
\]

This shows a decay of the fluctuations of the linear
statistics.\footnote{%
 The situation changes if we allow non-smooth test
 functions $ h $. In particular, for the flat model, the variance of
 the number of random zeroes in a smooth domain $ \Lambda $ is
 asymptotic to $ \sqrt{L} $ times the perimeter of $ \Lambda $, see
 \cite{FH99}. The difference between these 
 two types of behavior reflects high frequency oscillations of
 random zeroes which are not taken into account when the test
 function is smooth.}
Note that in the flat case formula \eqref{0.8L} was found in a paper
by Forrester and Honner \cite{FH99}, in which they also conjectured
the asymptotics in the elliptic case and checked it numerically. The
estimate $\Var Z_L(h) = O(1)$, $L\to\infty$, was obtained by Shiffman
and Zelditch in \cite{SZ99}.

\newtheorem*{MainT}{Main Theorem}
\begin{MainT}
In each of the three cases, the random variables
\[
\frac{ \sqrt L }{ \sqrt\kappa \| \Delta^* h \|_{L^2(m_*)} } \left(
Z_L (h) - L \frac1\pi \int_{\M} h \, dm^* \right)
\]
converge in distribution to $ \cN (0,1) $ for $ L \to \infty $.
\end{MainT}

The proof of Main Theorem is based on the asymptotic normality
theorem for non-linear functionals of Gaussian processes. Results
of this type  are known; usually, their proofs rely on the
classical method of moments combined with the diagram technique
(see Breuer and Major \cite{BM83}, and the references therein). In
Section~2 we shall prove another result in that spirit using a similar
strategy. In Section~3, we deduce Main Theorem from this result and
prove the asymptotics (\ref{0.8L}).

For real random polynomials $ P_N (x) = \sum_{k=0}^N \xi_k x^k $
where $\xi_k$ are independent identically distributed (real-valued)
random variables such that $ \Ex\xi_k = 0 $,
$ \Pr{ \xi_k=0 }=0 $, and $ \Ex |\xi_k|^{2+\epsilon} <
\infty $, Maslova \cite{Ma74a}, \cite{Ma74b} evaluated the variance of
the number of real zeroes of $ P_N $, and proved the corresponding
version of asymptotic normality. Her methods are
quite different from ours. Probably, Main Theorem persists for more
general models like zeroes of random holomorphic sections of high
powers of Hermitian line bundles over K\"ahler manifolds
extensively studied by Bleher, Shiffman and Zelditch \cite{SZ99,
BSZ00a, BSZ00b}.

\medskip\noindent\textsc{Three toy models for the flat CAZP.}
It is instructive to compare the flat CAZP with simpler (`toy') models
of random point processes in the plane, especially, random
perturbations of a lattice. \emph{The first toy model:} each point of
the lattice $ \sqrt{\frac\pi2} \Z^2 = \{ \sqrt{\frac\pi2} (k+li) : k,l
\in \Z \} $ is deleted at random, independently of others, with
probability $ 1/2 $; the remaining points are a random set. The
corresponding smooth linear statistics $ Z_L^{(1)} (h) = \int h
(z/\sqrt{L}) \, dn^{(1)} (z) $ (where $ n^{(1)} $ is the counting
measure on the random set) is asymptotically normal with parameters $
\Ex Z^{(1)}_L (h) = L \pi^{-1} \int h \, dm $ and $ \Var Z^{(1)}_L (h)
= \pi^{-1} \int h^2 (z/\sqrt{L}) \, dm(z) = L \pi^{-1} \| h
\|^2_{L^2(m)} $, which is quite different from (\ref{0.8L}): the
variance grows, not decays; and $ \| h \| $ appears, not $ \| \Delta h
\| $.

\emph{The second toy model:} points of the lattice $ \sqrt\pi \Z^2 $
move independently, forming the random set $ \{ \sqrt{\pi} (k+li) +
c \eta_{k,l} :\, k,l \in \Z \} $ where $ \eta_{k,l} $ are independent
standard complex Gaussian random variables, and $ c \in (0,\infty) $
is a parameter. The corresponding $ Z^{(2)}_L (h) = \int h
(z/\sqrt{L}) \, dn^{(2)} (z) $ is asymptotically normal with
parameters $ \Ex Z^{(2)}_L (h) = (1+o(1)) L \pi^{-1} \int h \, dm $
and $ \Var Z^{(2)}_L (h) = (c^2+o(1)) (2\pi)^{-1} \| \nabla h
\|^2_{L^2(m)} $ (where $ \nabla $ stands for gradient). This time, the
variance does not grow, but still, it does not decay; and $ \| \nabla
h \| $ is not $ \| \Delta h \| $.

\emph{The third toy model} reaches asymptotic similarity to the flat
CAZP, but is more complicated. Lattice points are initially aggregated
into three-point clusters, and each cluster scatters in a special
(equiangular and equidistant) way. Namely, we consider the random set
\[
\{ \sqrt{3\pi} (k+li) + c e^{2\pi im/3} \eta_{k,l} :\, k,l \in \Z;
\, m=0,1,2\},
\]
where $ \eta_{k,l} $ and $ c $ are as before. The
corresponding $ Z^{(3)}_L (h) = \int h (z/\sqrt{L}) \, dn^{(3)} (z) $
is asymptotically normal with parameters $ \Ex Z^{(3)}_L (h) =
(1+o(1)) L \pi^{-1} \int h \, dm $ and $ \Var Z^{(3)}_L (h) = (1+o(1))
L^{-1} 3 c^4 (16\pi)^{-1} \| \Delta h \|^2_{L^2(m)}$. This time, we
can mimic (\ref{0.8L}) by choosing $ c = 2 (\pi\kappa/3)^{1/4} $.

We come to a vague idea of two conservation laws for random zeroes:
mass conservation and center-of-mass conservation. Mass is conserved
in the second and third toy models; center-of-mass is conserved in the
third toy model only. Both should be conserved in the flat CAZP, in
some sense (to be understood). An attempt to understand the mass
conservation is made in our next work \cite{SoTsi2}. The center-of-mass
conservation remains unexplored.

About the three directions ($ e^{2\pi im/3} : m=0,1,2 $): four or more
directions may be used equally well, but two directions are not
enough. Indeed, every quadratic form $ Q : \C \to \R $ satisfies $
\frac1n \sum_{m=0}^{n-1} Q (e^{2\pi im/n}) = \frac14 \Delta Q(0) $
provided that $ n \ge 3 $; however, $ Q(1) + Q(-1) $ is not
proportional to $ \Delta Q(0) $.

Maybe, the failure of two-point clusters hints at a third conservation
law.

\medskip\noindent\textsc{Acknowledgment.}
We thank Leonid Pastur, Leonid Polterovich and Zeev Rudnik for useful
discussions.


%% file: sect1.tex
\section{Geometrical description of models}

By the geometrical description of CAZP we mean something like this:
\begin{quote}
The CAZP process on a complex 1-dimensional analytic manifold is the
intersection of its isometric embedding into a projective space with a
random hyperplane.
\end{quote}
However, we do not formalize the description; we only explain, how it
works in the three models considered. For related more advanced
theories, see the works of Gromov \cite[Sect.~3.3]{Gr90}, Shub and Smale
\cite{ShSm}, and Bleher, Shiffman and Zelditch \cite{BSZ01}.

\subsection{Three homogeneous spaces}

Each model is based on a
simply connected homogeneous space $ \M $; it may be thought of
as a real 2-dimensional manifold or a complex 1-dimensional analytic
manifold.\footnote{%
 The choice of a complex structure on our real manifold does not
 introduce arbitrariness; there exist only two invariant complex
 structures, conjugate to each other, both leading to the same CAZP.}
The former may be embedded into a Euclidean space, the
latter --- into the projective space $ P(\C^2) $ of
1-dimensional subspaces of $ \C^2 $. The complex plane $ \C $
is embedded into $ P(\C^2) $ by $ z \mapsto \{ (u,uz) : u \in \C \} $;
in order to cover the whole $ P(\C^2) $, one additional point is
needed, $ \infty \mapsto \{ (0,u) : u \in \C \} $. The symmetry group
$ G $ in the real case is a subgroup of the group of motions of the
Euclidean space. In the complex case, $ G $ is a subgroup of the group
of projective transformations of $ P(\C^2) $. The latter is covered by
the group $ SL(2) $, and we may take $ G \subset SL(2) $; such an
action need not be effective (since $ (-1) $ acts trivially), which is
harmless.

\medskip\noindent\textsc{Elliptic model.}
The real manifold is the sphere
\[
\S^2 = \{ (x_0,x_1,x_2) : x_0^2 + x_1^2 + x_2^2 = 1 \} \subset \R^3,
\]
the symmetry group being $ SO(3)
$ (orientation-preserving rotations of $ \R^3 $). The complex manifold
is the whole $ P(\C^2) $, with $ G = SU(2) $; or alternatively, $ \C
\cup \{\infty\} $ with transformations $ z \mapsto (az+b) / (cz+d) $,
$ \( \begin{smallmatrix} a & b \\ c & d \end{smallmatrix} \) \in SU(2)
$, preserving the spherical metric $ |dz| / (1+|z|^2) $. The
correspondence (well-known as the stereographic projection) is $ z =
(x_1 + ix_2) / (1-x_0) $. See also \cite[Sect.~2]{Le99} and
\cite[Chap.~4]{Pe}. 

\medskip\noindent\textsc{Flat model.}
The real manifold is the plane $ \R^2 $ with the group of
(orientation-preserving) motions, that is, shifts and rotations. The
complex manifold is the complex plane $ \C $ with transformations $ z
\mapsto e^{i\phi} (z+u) $, preserving the Euclidean metric $ |dz| $;
or alternatively, $ P(\C^2) $ minus a single point, with
transformations preserving the point. The correspondence is just $ z =
x_1 + i x_2 $.

\medskip\noindent\textsc{Hyperbolic model.}
The real manifold is $ \{ (x_0,x_1,x_2) : x_0 = \sqrt{ 1 + x_1^2 +
x_2^2 } \} \linebreak[0]
\subset \R^3 $, the upper sheet of a hyperboloid (in other
words, a pseudosphere), the symmetry group being the connected
component of $ SO(2,1) $ (rotations of $ \R^3 $ that preserve the
sheet and orientation). The complex manifold is the disc $ \D = \{ z
\in \C : |z| < 1 \} $, with transformations $ z \mapsto (az+b) /
(\overline b z + \overline a) $, $ |a|^2 - |b|^2 = 1 $, preserving the
hyperbolic metric $ |dz| / (1-|z|^2) $; or alternatively, the
corresponding part of $ P(\C^2) $, with $ G = SU(1,1) $ (isomorphic to
$ SL(2,\R) $). The correspondence is $ z = (x_1 + ix_2) / (1+x_0)
$. See also \cite[Sect.~2]{Le99} and \cite[Chap.~5]{Pe}. 

\subsection{Enlarging symmetry}

The symmetry group $ G $ acts transitively on $ \M $. It acts also on
discrete subsets of $ \M $, but not transitively. A $ G $-invariant
random discrete subset of $ \M $ is far from being unique in
distribution.

The key ingredient of the construction is an embedding $ \iota $ of $ \M $
into a high-dimensional projective space $ P(\C^n) $, $ n \in \{
1,2,\dots \} \cup \{\infty\} $, and use of the high-dimensional
symmetry group $ U(n) $ of rotations of $ \C^n $. Of course, the image
$ \iota (\M) $ is not $ U(n) $-invariant. However, we use $ U(n)
$-invariance for determining a probability measure on
hyperplanes. Choosing at random a hyperplane, we observe its
intersection with $ \iota (\M) $. The case $ n=\infty $ does involve
some technicalities; we will return to the point later.

The embedding must be $ ( G, U(n) ) $-invariant in the following
sense: for every $ g \in G $ there exists $ U \in U(n) $ such that the
diagram
\begin{equation}\label{1.1}
\begin{gathered}
\xymatrix{
 \M \ar[r]^g \ar[d]_\iota & \M \ar[d]_\iota \\
 P(\C^n) \ar[r]^{P(U)} & P(\C^n)
}
\end{gathered}
\end{equation}
is commutative. Surprisingly, such an embedding is unique (up to a
rotation)! Uniqueness of $ \iota $ leads to a model, (statistically)
unique in the same sense as the Poisson process (recall the quote from
Hannay in the Introduction).

Now we switch from a heuristic to rigorous style.

Let $ H $ be a Hilbert space over $ \C $, either finite-dimensional,
or infinite-dimensional and separable. The projective space $ P(H) $
is, by definition, the set of all one-dimensional subspaces of $ H
$. Each non-zero vector $ x \in H \setminus \{0\} $ spans such a
subspace $ P(x) \in P(H) $. A transformation $ P(U) : P(H) \to P(H) $
corresponds to every unitary operator $ U \in U(H) $; namely, $ P(U)
\( P(x) \) = P(Ux) $. Note that $ P(e^{i\phi}U) = P(U) $. Two maps $
\iota, \iota' : \M \to P(H) $ are called unitarily equivalent, if $ \iota' 
= P(U)
\circ \iota $ for some $ U \in U(H) $. A map $ \iota : \M \to P(H) $ is called
holomorphic, if locally (in a neighborhood of any point of $ \M $) it
is the composition of some holomorphic map $ \M \to H
\setminus \{0\} $ and the canonical projection $ H \setminus \{0\} \to
P(H) $, $ x \mapsto P(x) $.

The well-known Fubini-Study metric\footnote{%
 Or rather, the geodesic metric corresponding to the Fubini-Study
 metric tensor, up to a coefficient.}
on $ P(H) $,
\[
\dist \( P(x), P(y) \) = \arccos \frac{ | \ip x y | }{ \|x\| \|y\| }
 \, ,
\]
is $ U(n) $-invariant. Given a one-to-one map $ \iota : \M \to P(H) $, we 
get
a metric on $ \M $,
\[
\dist_\iota (z,z') = \dist \( \iota (z), \iota (z') \) \, ;
\]
assuming smoothness of $ \iota $ and taking $ z' $ infinitesimally close
to $ z $ we get a tensor field $ g_\iota $ on $ \M $, the Riemannian metric
induced by $ \iota $. If $ \iota $ and $ \iota' $ are unitarily
equivalent, then they induce the same Riemannian metric (since $
\dist_\iota = \dist_{\iota'} $).

\newtheorem*{CalabiT}{Calabi's rigidity theorem}
\begin{CalabiT}
If two holomorphic embeddings\footnote{%
 `Embedding' means one-to-one. The theorem holds also for immersions,
 but we need only embeddings.}
of a complex manifold into $ P(H) $ induce the same Riemannian metric on
$ \M $, then they are unitarily equivalent. \textup{(\cite[Th.~9]{Ca53}, see
 also \cite[Sect.~4]{SV}.)}
\end{CalabiT}

If a holomorphic embedding $ \iota : \M \to P(H) $ is $ ( G, U(n) )
$-invariant (recall \eqref{1.1}), then it induces a $ G $-invariant
Riemannian metric on $ \M $. However, such a metric is unique up to a
coefficient (since $ G $ can rotate $ \M $ around any given point). It
means that \emph{the coefficient is the only possible parameter of a
holomorphic $ ( G, U(n) ) $-invariant embedding $ \iota : \M \to P(H) $}
(treated up to unitary equivalence).

Another implication of Calabi's theorem: \emph{if a holomorphic
embedding $ \iota : \M \to P(H) $ induces a $ G $-invariant Riemannian 
metric on $ \M $, then $ \iota $ is $ ( G, U(n) ) $-invariant.}

By the way, $ ( G, U(n) ) $-invariance leads to a \emph{projective}
representation of $ G $ (that is, a homomorphism from $ G $ to the
factor group of $ U(H) $ by its center), not necessarily a homomorphism
$ G \to U(H) $.

\subsection{Relation to the analytical description}

In order to finish the description, it remains
to write down the corresponding embeddings $ \iota $ in each of the
three cases, to check invariance of induced metrics, and to explain
the rotation invariant choice of the random hyperplane in the case
when the dimension of $ H $ is infinite.

All invariant metrics on $ \M $ are proportional to the one
mentoned in the Introduction. The only freedom left to us is to choose the
numerical coefficient $ \sqrt L $ of the invariant metric on $ \M $.

\medskip\noindent\textsc{Elliptic case.}
For the sphere $ \S^2 $, the existence of such an embedding depends on
the parameter $ L $. Only for $ L = 1,2,3,\dots $ the invariant metric
on $ \S^2 $ is embeddable into $ P(H) $.\footnote{
 One can easily get this using the Edelman-Kostlan formula: if $
 L $ is non-integer, then we would get a non-integer answer for
 the average number of zeroes of our random function on $ \S^2 $.}
For such $ L $, an embedding $ \S^2 \to
P(\C^{L+1}) $ is well-known to physicists as the system of
spin-$ J $ coherent states ($ J = L/2 $), see \cite{Ha96},
\cite[Sect.~4.3]{Pe}, \cite{SV}. Sometimes, mathematicians call the
wave function of a coherent state ``the Szeg\"o kernel''
(see \cite{BSZ00b}). Treated up to rotations of $ \C^{L+1} $, the
map is unique. By rigidity, every such embedding uses only a
finite-dimensional subspace of $ H $.

The embedding $ \S^2 \to P(\C^{L+1}) $ may be described as the
composition $ \S^2 \to \C \to \C^{L+1} \setminus \{0\} \to P(\C^{L+1})
$, where $ \S^2 \to \C $ is the stereographic projection mentioned in
Sect.~1.1, $ \C^{L+1} \setminus \{0\} \to P(\C^{L+1}) $ is the
canonical projection $ x \mapsto P(x) $, and $ \ti \iota : \C \to \C^{L+1}
\setminus \{0\} $ is given by
\[
\ti \iota (z) = \Big( {\textstyle \sqrt{\binom L 0}, \sqrt{\binom L 1}
\, z, \dots, \sqrt{\binom L L} \, z^L } \Big) \, ,
\]
which evidently corresponds to \eqref{elliptic}.

The induced metric is easy to calculate:
\begin{gather*}
\ip{ \ti \iota (z) }{ \ti \iota (z') } = \sum_{k=0}^L { \textstyle
\sqrt{\binom L k} \, z^k
 \sqrt{\binom L k} \, \bar z^{\prime k} } = (1+z \bar z')^L \, ; \\
\dist_\iota (z,z') = \arccos \frac{ | \ip{ \ti \iota (z) }{ \ti \iota (z') } | }{
\| \ti \iota (z) \| \| \ti \iota (z') \| } = \arccos \bigg( \frac{ | 1+z \bar
z' | }{ \sqrt{1+|z|^2} \sqrt{1+|z'|^2} } \bigg)^L \, .
\end{gather*}
$ G $-invariance of the metric for every $ L $ follows immediately
from its $ G $-invariance for $ L=1 $; however, for $ L=1 $ the
embedding $ \C \to P(\C^{L+1}) $, given by $ z \mapsto \{ (u,uz) : u
\in \C \} $, is already familiar to us (recall the beginning of
Sect.~1.1). An explicit calculation gives for $ \De z \to 0 $
\begin{multline*}
\dist_\iota \bigg( z - \frac12 \De z, z + \frac12 \De z \bigg) \\
= \arccos \bigg( 1 - \frac{ |\De z|^2 }{ (1+|z|^2+\frac14 |\De z|^2)^2
 - \Re^2 (\bar z \De z) } \bigg)^{L/2} \\
= \arccos \bigg( 1 - \frac L 2 \frac{ |\De z|^2 }{ (1+|z|^2)^2 }
 (1+o(1)) \bigg) = \sqrt L \frac{ |\De z| }{ 1+|z|^2 } (1+o(1)) \, .
\end{multline*}

\medskip\noindent\textsc{Flat case.}
If the needed embedding $ \iota : \C \to P(H) $ exists for $ L=1 $, then
for every $ L \in (0,\infty) $ the map $ z \mapsto \iota (L^{1/2}z) $ fits;
we restrict ourselves to $ L=1 $. The construction of $ \iota $ is
well-known to physicists as the usual system of coherent states
(Schr\"odinger, von Neumann, Klauder, et
al.), see \cite[Chapter 1]{Pe}. Its explicit form is given
(similarly to the elliptic case) by the composition $ \C \to H
\setminus \{0\} \to P(H) $ where the first map $ \ti \iota : \C \to H
\setminus \{0\} $ is given by
\[
\ti \iota (z) = \bigg( 1, z, \frac{z^2}{\sqrt{2!}}, \frac{z^3}{\sqrt{3!}},
\dots \bigg) \in l^2 = H \, ,
\]
cf. (\ref{flat}).

We have
\[
\ip{ \ti \iota (z) }{ \ti \iota (z') } = \sum_{k=0}^\infty
 \frac{z^k}{\sqrt{k!}} \frac{\bar z^{\prime k}}{\sqrt{k!}} = \exp ( z
 \bar z' ) \, ;
\]
\begin{multline*}
\dist_\iota (z,z') = \arccos \exp \Big( \Re (z\bar z') - \frac12 |z|^2 -
 \frac12 |z'|^2 \Big) \\
= \arccos \exp \Big( -\frac12 |z-z'|^2 \Big) \, ;
\end{multline*}
\[
\dist_\iota ( z, z + \De z ) = \arccos \bigg( 1 - \frac12 |\De z|^2 (1+o(1))
\bigg) = |\De z| (1+o(1)) \, .
\]

\medskip\noindent\textsc{Hyperbolic case.}
The corresponding system of coherent states is usually investigated in
terms of unitary representations of the group $ G = SU(1,1) $, see
\cite[Chap.~5]{Pe}, \cite[Sect.~2]{Le99}, which leads to a special
treatment of integer values of $ L $. In our approach,
\emph{projective} representations appear irrespective of unitary
representations, and $ L $ runs over $ (0,\infty) $.

As before, $ \iota (z) = P(\ti \iota(z)) $, but now $ z \in \D $, and
$ \ti \iota : \D \to l^2 $ is given by
\[
\ti \iota (z) = \Big( {\textstyle 1, \sqrt{L} \, z,
\sqrt{\frac{L(L+1)}{2!}} \, z^2, \sqrt{\frac{L(L+1)(L+2)}{3!}} \, z^3,
\dots } \Big) \, ,
\]
cf.~(\ref{hyperbolic}). The corresponding embedding of the
hyperbolic plane to $ P(H) $ was known already to Bieberbach in
1932. We have
\begin{gather*}
\ip{ \ti \iota (z) }{ \ti \iota (z') } = \sum_{k=0}^\infty {\textstyle
\sqrt{\binom{L+k}k} \, z^k \sqrt{\binom{L+k}k} \, \bar z^{\prime k} }
= (1-z \bar z')^{-L} \, ;
 \\
\dist_\iota (z,z') = \arccos \bigg( \frac{ | 1-z \bar z' | }{
 \sqrt{1-|z|^2} \sqrt{1-|z'|^2} } \bigg)^{-L} \, ;
\end{gather*}
\vspace{-4mm}
\begin{multline*}
\dist_\iota \bigg( z - \frac12 \De z, \, z + \frac12 \De z \bigg) \\
= \arccos
 \bigg( 1 + \frac{ |\De z|^2 }{ (1-|z|^2-\frac14 |\De z|^2)^2 - \Re^2
 (\bar z \De z) } \bigg)^{-L/2} \\
= \arccos \bigg( 1 - \frac L 2 \frac{ |\De z|^2 }{ (1-|z|^2)^2 }
 (1+o(1)) \bigg) = \sqrt L \frac{ |\De z| }{ 1-|z|^2 } (1+o(1)) \, .
\end{multline*}

\subsection{Random hyperplane in Hilbert space?}

In $ \C^n $, the only $ U(n) $-invariant distribution on hyperplanes,
the uniform distribution, may be represented via the normal (to the
hyperplane) vector $ (\zeta_1,\dots,\zeta_n) $ distributed uniformly
on the sphere. The normal distribution for $ (\zeta_1,\dots,\zeta_n) $
can be used instead (as well as any spherically invariant
distribution).

In a Hilbert space $ H = l^2 $ there is no rotation-invariant
probability measure on hyperplanes (nor Lebesgue measure in $ H
$). Nevertheless, an infinite sequence of independent normal $ \cN_\C
(0,1) $ random variables $ \zeta_k $ can be used. Of course, it is not
a random element of $ H $, since the event $ \sum |\zeta_k|^2 < \infty
$ is of probability $ 0 $. However, for every $ x = (c_1, c_2, \dots)
\in l^2 $ the series $ \sum c_k \bar\zeta_k $ converges a.s., and we
may \emph{denote} its sum by $ \ip x \zeta $. The `bad' set of zero
probability, on which the series does not converge, depends on $ x $;
the union over all $ x \in H $ is not a set of zero
probability.\footnote{%
 In fact, it is the event $ \sum |\zeta_k|^2 = \infty $, of
 probability $ 1 $.}
We cannot choose $ \zeta $ at random and speak about `the function $ x
\mapsto \ip x \zeta $ on $ H $'.

What we really need, is the function $ z \mapsto \ip{ \ti \iota (z) } \zeta
$ for $ z \in \M $, that is, the function $ x \mapsto \ip x \zeta $
for $ x \in \ti \iota (\M) $ only. The corresponding series converges
a.s. for all these $ x $ simultaneously. Thus, `the random hyperplane
of $ H $' is ill-defined, but still, its intersection with $ \ti \iota
(\M)
$ is well-defined.

It remains to explain unitary invariance of the construction described
above. Random variables $ \ip x \zeta $ can be used simultaneously for
a \emph{countable} set of points $ x $. In particular, for every
orthonormal basis $ (e_1,e_2,\dots) $ of $ H $, the sequence of random
variables $ \ip{e_k}\zeta $ is well-defined. Some reflection shows that
the joint distribution of these $ \ip{e_k}\zeta $ does not differ from
that of $ \zeta_k $. It follows easily that our construction can start
with an arbitrary basis, resulting in the same distribution of the
random function $ z \mapsto \ip{ \ti \iota (z) } \zeta $. See also
\cite[Sect.~1.3]{Ja}.


%% file: sect2.tex
\section{Asymptotic normality for non-linear functionals of Gaussian
 processes}

\subsection{The result}

Let $T$ be a measure space endowed with a finite
positive measure $\mu$, $\mu (T)=1$. A complex-valued Gaussian process
on $T$ may be defined as
\begin{equation}\label{eq1.1}
w(\om, t) = \sum_k \zeta_k (\om) f_k(t)
\end{equation}
where $f_k : T\to \C$ are measurable functions such that
\[
\sum_k |f_k(t)|^2 <\infty \qquad \text{for all } t\in T,
\]
and $\zeta_k = \xi_k +i\eta_k$ are independent standard complex
Gaussian variables; i.e. $\zeta_k \sim  \cN_\C(0,1) $. The latter means
that $\xi_k$ and $\eta_k$ are
independent centered Gaussian (real) variables with variance
$\frac12$; then $\Ex\zeta_k = 0$ and $\Ex |\zeta_k|^2
= 1$. We restrict ourselves to the case
\[
\sum_k |f_k(t)|^2 = 1 \qquad \text{for all } t\in T\,.
\]
Then $w(t)\sim \cN_{\C}(0,1)$ for all $t\in T$.

Now, a few words about the convergence of the series \eqref{eq1.1}. Treating
each term
$\zeta_k(\om)f_k(t)$ as an element of the space $\cL^2 =
L^2\left( (\Om, P) \times (T, \mu) \right)$, we have
\[
\|\zeta_k f_k\|_{\cL^2} = \|\zeta_k\|_{L^2(\Om, P)}
\|f_k\|_{L^2(T, \mu)} = \|f_k\|_{L^2(T, \mu)}
\]
and the terms are pairwise orthogonal. Therefore, the series \eqref{eq1.1}
converges in the space\footnote{%
 Moreover, for each $t\in T$, the series
 \eqref{eq1.1} converges in $L^2(\Om, P)$ (for an obvious reason), and
 by the Kolmogorov-Khinchin theorem \cite[Ch.~3]{Ka85} the series converges
 in $L^2(T, \mu)$ for almost all $\om\in \Om$.}
$\cL^2$.

In what follows, we always assume that the sum of the series \eqref{eq1.1}
is not the zero function in the space $\cL^2$.

The correlation function $\rho: T\times T \to \C$ of the process
$w(t)$ equals
\[
\rho (s,t) \stackrel{\mathrm{def}}{=} \Ex \left\{ w(s) \overline{w(t)}
\right\} = \sum_k f_k(s) \overline{f_k(t)}\,.
\]
Clearly, $|\rho (s,t)|\le 1$ and $\rho (t,t)=1$.

Consider a sequence of complex Gaussian processes $w_n$ with the
correlation functions $\rho_n (s,t)$ and define a sequence of random
variables
\[
Z_n = \int_T \phi (|w_n(t)|) \Theta (t) \, d\mu (t)
\]
where $\phi : \R_+\to \R$ is a measurable
function such that
\[
\int_0^\infty \phi^2(r)e^{-r^2/2}r\, dr <\infty,
\]
and $\Theta : T\to \R$ is a measurable bounded function.
We shall prove that under some natural assumptions on the
off-diagonal decay of the correlation functions $\rho_n (s,t)$
when $n\to\infty$, the random variables $Z_n$ are asymptotically
normal.

\begin{theorem}\label{2.2}
Suppose that for each $\alpha \in \N$
\begin{equation}\label{eq1.3}
\liminf_{n\to\infty} \frac{\iint_{T^2}
|\rho_n(s,t)|^{2\alpha} \Theta(s)\Theta(t)\,
d\mu(s)d\mu(t)}{\sup_{s\in T} \int_T |\rho_n(s,t)|\, d\mu(t)} >0,
\end{equation}
and that
\begin{equation}\label{eq1.4}
\lim_{n\to\infty} \sup_{s\in T} \int_T |\rho_n(s,t)|\, d\mu(t) = 0.
\end{equation}
Then the distributions of the random variables
\[
\frac{Z_n - \Ex Z_n}{\sqrt{\Var Z_n}}
\]
converge weakly to $\cN(0,1)$ for $n\to\infty$.

If $\phi$ is an increasing function, then it suffices to assume that
condition \eqref{eq1.3} holds only for $\alpha=1$.
\end{theorem}

\medskip\noindent{\bf Remarks:}

\smallskip\noindent{\bf (i) } The function $|\rho (s,t)|^{2\alpha}$
which appears in the integrand on the left-hand side of condition
\eqref{eq1.3} is a positive definite function on $T\times T$ (see formula
\eqref{eq1.6} below). Hence, for $\alpha\in\N$, the integral
$\iint_{T^2} |\rho (s,t)|^{2\alpha} \Theta (s) \Theta (t)\,
d\mu(s) d\mu (t)$ is always non-negative.

\smallskip
\begin{sloppypar}\noindent{\bf (ii) }
The role of condition \eqref{eq1.4} is to guarantee that
\[
\lim_{n\to\infty} \Var Z_n = 0.
\]
In fact, under assumption \eqref{eq1.3}, the sequences
\[
n\mapsto \mbox{Var} Z_n \quad \text{and} \quad
 n\mapsto \int_T |\rho_n(s,t)|\, d\mu(t)
\]
have the same decay (see expression \eqref{eq1.7} below).
\end{sloppypar}

\subsection{Some preliminaries}

We shall deal with the space $L^2_{\C}
\stackrel{\mathrm{def}}= L^2\left( \C, \cN_{\C}(0,1)
\right)$ and its subspaces $\mathcal P_m$ which consist of
the polynomials in $\zeta$ and $\overline\zeta$ of degrees
$\alpha$ and $\beta$ correspondingly, $\alpha+\beta=m$. The space
$L^2_\C$ has a polynomial basis
\[
\left\{ \frac1{\sqrt{\alpha! \beta!}} \Wi{ \zeta^\alpha
\overline\zeta^\beta } \right\}_{\alpha, \beta \in \Z_+}
\]
where $ \Wi{\zeta^\alpha \overline\zeta^\beta} $ is the
orthogonal projection of the polynomial $\zeta^\alpha
\overline\zeta^\beta$ onto the subspace $H^{:m:} \stackrel{\mathrm{def}}=
\mathcal P_m \ominus \mathcal P_{m-1}$, $m=\alpha+\beta$ \cite[Example
3.~32]{Ja}.\footnote{%
 The expression $ \Wi{ \zeta^\alpha
 \overline\zeta^\beta } $ is called the \emph{Wick product}.}
Thus any square
integrable function $\Phi$ of a random variable $\zeta\sim
\cN_\C(0,1)$ is of the form
\[
\Phi (\zeta) = \sum_{\alpha, \beta\in \Z_+} \frac{c_{\alpha
\beta} }{\sqrt{\alpha! \beta!}}
\Wi{ \zeta^\alpha\overline\zeta^\beta } \, , \qquad \|\Phi\|^2 =
\sum_{\alpha, \beta} |c_{\alpha \beta}|^2\,.
\]

In what follows, we shall deal only with the radial functions $\phi
(|\zeta|)$. In this case,
\begin{equation}\label{eq1.5}
\phi (|\zeta|) = \sum_{\alpha=0}^\infty
\frac{c_{2\alpha}}{\alpha!} \Wi{ |\zeta|^{2\alpha} }
\end{equation}
where $c_{2\alpha} \in\R$ and $\sum_\alpha c_{2\alpha}^2 =
\|\phi\|^2$. Indeed, the group of rotations $\zeta\mapsto
e^{i\theta} \zeta$ acts on the space $L^2_\C$ leaving
invariant the subspaces $H^{:m:}$. If $m$ is odd, then $H^{:m:}$
contains no radial functions; if $m=2\alpha$ is even, then the
subspace of radial polynomials in $H^{:m:}$ is one-dimensional and
is spanned by $ { \Wi{ |\zeta|^{2\alpha} } } =
{ \Wi{ \zeta^\alpha \overline\zeta^\alpha } }  $.

In the proof of Theorem~\ref{2.2}, we may assume without loss of generality that
\[
\int_{\C} \phi(|\zeta|) e^{-|\zeta|^2/2}\, dm(\zeta) =
0,
\]
that is, $\Ex Z_n = 0$. Then expansion \eqref{eq1.5} starts with $\alpha  =
\alpha_0 \ge 1$, $c_{2\alpha_0}\ne 0$.

Since
\begin{equation}\label{eq1.6}
\Ex \left( \frac1{\alpha!} \Wi{ |w(s)|^{2\alpha} }
\right) \left(\frac1{\beta!} \Wi{ |w(t)|^{2\beta} } \right) =
\begin{cases}
|\rho(s,t)|^{2\alpha} & \text{if $\alpha = \beta$}, \\
0 & \text{otherwise}
\end{cases}
\end{equation}
\cite[Theorem~3.9]{Ja} (see also Example at the end of the next
subsection), we have
\[
\Ex \left\{ \phi(|w_n(s)|) \phi(|w_n(t)|)\right\} =
\sum_{\alpha\ge \alpha_0} c_{2\alpha}^2 |\rho_n(s,t)|^{2\alpha} \, ,
\]
and
\begin{multline}\label{eq1.7}
\Ex Z_n^2 = \Ex \left( \int_T \phi (|w_n(t)|)
 \Theta (t) \, d\mu (t) \right)^2 \\
= \iint_{T^2} \Ex \left\{ \phi(|w_n(s)|)
 \phi(|w_n(t)|)\right\} \Theta(s) \Theta(t) \, d\mu(s) d\mu(t) \\
= \sum_{\alpha\ge \alpha_0} c_{2\alpha}^2 \iint_{T^2}
 |\rho_n(s,t)|^{2\alpha} \Theta(s) \Theta(t) \, d\mu(s) d\mu(t)
 \, .
\end{multline}

The proof of Theorem~\ref{2.2} uses the classical method of moments,
though it cannot be applied directly to the random variables
$Z_n$ since $\phi$ (and therefore, $Z_n$) need not have more
than two moments. First, we approximate $\phi$ by the
polynomials
\[
\phi_m (|\zeta|) = \sum_{\alpha\le m}
\frac{c_{2\alpha}}{\alpha!} \Wi{ |\zeta|^{2\alpha} } \, \qquad
m\ge \alpha_0,
\]
and $Z_n$ by the random variables
\[
Z_{n,m} = \int_T \phi_m(|w_n(t)|) \Theta(t) \, d\mu(t)\,.
\]
Then the moment method will be applied to $Z_{n,m}$.

Applying formula \eqref{eq1.7} to $Z_n-Z_{n,m}$ (that is, replacing the
function $\phi$ by $\phi-\phi_m$), we get
\begin{multline*}
\Ex \left( Z_n-Z_{n,m}\right)^2 = \sum_{\alpha\ge m+1}
 c_{2\alpha}^2 \iint_{T^2} |\rho_n(s,t)|^{2\alpha}
 \Theta(s) \Theta(t)
 \, d\mu (s) d\mu (t) \\
\le \|\Theta\|_{L^\infty (T, \mu)}^2 \sum_{\alpha\ge m+1}
 c_{2\alpha}^2 \sup_{s\in T}\int_T |\rho_n(s,t)|
 \, d\mu (t) \\
\stackrel{\eqref{eq1.3}}\le \frac{\|\Theta\|_{L^\infty (T,
\mu)}^2}{\eta(\alpha_0)} \, \|\phi - \phi_m\|^2_{L^2_\C} \,
 \iint_{T^2} |\rho_n(s,t)|^{2\alpha_0} \Theta (s) \Theta
(t) \, d\mu (s) d\mu (t)
\end{multline*}
where $\eta (\alpha)$ denotes the expression on the left-hand side of
\eqref{eq1.3}. Using again formula \eqref{eq1.7}, we finally get
\[
\Ex \left(Z_n-Z_{n,m}\right)^2 \le
\frac{\|\Theta\|^2_{L^\infty (T, \mu)}}{c_{2\alpha_0}^2 \eta
(\alpha_0)} \, \|\phi - \phi_m\|^2_{L^2_\C} \,
\Ex Z_n^2\,.
\]
Hence
\[
\lim_{m\to\infty} \sup_n \Ex \left(
\frac{Z_n}{\sqrt{\hbox{Var} Z_n}} -
\frac{Z_{n,m}}{\sqrt{\hbox{Var} Z_{n,m}}} \right)^2 = 0\,,
\]
and if for each fixed $m\ge \alpha_0$ the random variables
$ Z_{n,m} / \sqrt{\Var Z_{n,m}} $ converge in
distribution to $\cN(0,1)$ when $n\to\infty$, then the
random variables $ Z_n / \sqrt{\Var Z_n} $ have the same
property.

From now on, we always assume that $\phi$ is a polynomial (and
hence $Z_n$ has the moments of all orders).

\subsection{More preliminaries (the diagram formula)}

The next step is to evaluate the moments
$ \Ex Z_n^p $, $p\in\N$, and to compare them with
$\left( \Ex Z_n^2\right)^{p/2}$ times the moments of the
standard normal distribution.

We have
\begin{multline}\label{eq1.8}
\Ex Z_n^p = \idotsint_{T^p} \Ex \left\{
 \prod_{j=1}^p \phi (|w_n(t_j)|) \right\} \prod_{j=1}^p \Theta
 (t_j) d\mu (t_j) \\
\stackrel{\eqref{eq1.5}}= \sum_{\alpha_1,\dots,\alpha_p\ge \alpha_0}
\frac{c_{2\alpha_1}\dots c_{2\alpha_p}}{\alpha_1!\dots \alpha_p!}
\idotsint_{T^p} \Ex \left\{ \prod_{j=1}^p
\Wi{ |w_n(t_j)|^{2\alpha_j} } \right\} \prod_{j=1}^p \Theta
(t_j) d\mu (t_j) \,.
\end{multline}
We compute the integrand  $\Ex \left\{ \prod_{j=1}^p
\Wi{ |w_n(t_j)|^{2\alpha_j} } \right\}$ using the diagram
technique \cite[Chapter 3]{Ja}.

Fix the exponents $\alpha_1, \dots , \alpha_p$. A {\em diagram}
$\gamma$ is a graph with $ 2(\alpha_1+\dots +\alpha_p) $ vertices
labeled by the indices $1, \overline 1, 2, \overline 2, \dots , p,
\overline p$ ($\alpha_j$ vertices are labeled by $j$ and other
$\alpha_j$ vertices are labeled by $\overline j$), and each
vertex has degree one (i.e., the edges have no common end points).
The edges may connect only the vertices labeled by $i$ and
$\overline j$ with $i\ne j$. The set of all such graphs is denoted
by $\Gamma (\alpha_1, \dots , \alpha_p)$. For some choices of
$\alpha_1, \dots , \alpha_p$, $\Gamma(\alpha_1, \dots , \alpha_p)$ may
be the empty set, for example, $\Gamma (\alpha_1,
\alpha_2)=\emptyset$ iff $\alpha_1\ne \alpha_2$, and $\Gamma
(\alpha_1, \alpha_2, \alpha_3) = \emptyset$ if $ \alpha_1 >
\alpha_2 + \alpha_3 $. The {\em value} of the diagram $\gamma$
equals
\[
V_\gamma (t_1, \dots , t_p) \stackrel{\mathrm{def}}= \prod_{(i,\overline j)\in
e(\gamma)} \rho (t_i, t_j)
\]
where the product is taken over all edges of $\gamma$. In this
notation, the diagram formula \cite[Theorem~3.12]{Ja}
reads\footnote{%
 The diagrams we consider do not contain edges
 which join the vertices labeled by $i$ and $j$ (and by $\overline i$
 and $\overline j$) since $\Ex\left\{ w(t_i)w(t_j)\right\} =
 \Ex\left\{ \overline{w(t_i)}\, \overline{w(t_j)}\right\} = 0
 $.}
\begin{equation}\label{eq1.9}
\Ex
\left\{ \prod_{j=1}^p \Wi{ |w_n(t_j)|^{2\alpha_j} } \right\}
= \sum_{\gamma \in \Gamma (\alpha_1, \dots , \alpha_p)} V_\gamma (t_1, \dots ,
t_p)
\end{equation}
(as usual, summation over the empty set means that the sum has zero
value).

\begin{example}
Consider the simplest case $p=2$. If
$\alpha_1\ne \alpha_2$, then $\Gamma(\alpha_1, \alpha_2)= \emptyset$, and
$ \Ex \left\{ \Wi{ |w(t_1)|^{2\alpha_1} }
\Wi{ |w(t_2)|^{2\alpha_2} } \right\} = 0 $. Now, suppose
$ \alpha_1 = \alpha_2 = \alpha $. Then we can glue together the vertices
labeled by the same indices. We get a graph with four vertices and two
edges of multiplicity $ \alpha $:
\begin{equation}\label{fig1}
\begin{gathered}\includegraphics{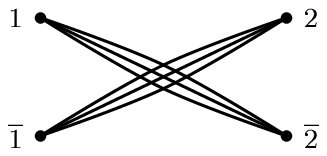}\end{gathered}
\end{equation}
The edges connecting $1$ and $\overline 2$ contribute the factor
$\rho(t_1, t_2)^\alpha$, the other $\alpha$ edges  contribute the
factor $\rho (t_2, t_1)^\alpha = \overline{ \rho (t_1,
t_2)}\,^\alpha$. Thus the value of the diagram is $|\rho (t_1,
t_2)|^{2\alpha}$, and $ \Ex \left\{
\Wi{ |w(t_1)|^{2\alpha} } \Wi{ |w(t_2)|^{2\alpha} } \right\}
= \sharp \Gamma (\alpha, \alpha) \cdot |\rho (t_1, t_2)|^{2\alpha}
$. It remains to find the total number of diagrams in $\Gamma
(\alpha, \alpha)$.

All the diagrams in $\Gamma (\alpha, \alpha)$ can be obtained from
the fixed one by permutation of $\alpha$ vertices labeled by $1$,
and by another independent permutation of $\alpha$ vertices
labeled by $2$. Therefore, $\sharp \Gamma (\alpha, \alpha) =
(\alpha!)^2$, and we recover formula \eqref{eq1.6}.
\end{example}

\subsection{The main argument}

For $p$ even, there are diagrams with a simple
structure whose total contribution to $\Ex Z_n^p$ {\em equals}
$\left( \Ex Z_n^2  \right)^{p/2} \Ex \xi^p $, where $\xi \sim
\cN(0,1)$. Contribution of the other diagrams is negligible.

\begin{definition}
A diagram $\gamma $ is called
{\em regular} if the set $\{1, 2, \dots , p\}$ is split into $q=p/2$
pairs and there are no edges between the vertices $i$ and
$\overline j$ if $i$ and $j$ belong to different pairs. Otherwise,
the diagram is called {\em irregular}.
\end{definition}

In other words, the diagram is regular if after glueing together
the vertices labeled by the same index it becomes a disjoint
union of $q=p/2$ ``elementary diagrams'' drawn on \eqref{fig1} (having,
generally speaking, different multiplicities of the edges).
\[
\begin{gathered}\includegraphics{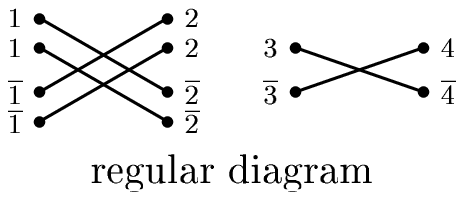}\quad
                \includegraphics{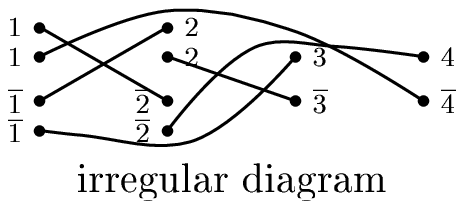}\end{gathered}
\]

Suppose $p$ is even and the diagram $\gamma $ is regular. For each
$j$, we glue together all the vertices labeled by $j$ and
$\overline j$, obtaining a {\em reduced diagram} $\gamma^*$; i.e.,
a graph with $p$ vertices and multiple edges which is split into
$q=p/2$ connected components.

For example, if $\gamma $ is related to the partition
\begin{equation}\label{eq1.10}
\{1, 2, \dots , p \} = \{1, 2\} \cup \{3, 4\} \cup\dots  \cup \{p-1,
p\} \, ,
\end{equation}
then the reduced diagram $\gamma^*$ is
\[
\begin{gathered}\includegraphics{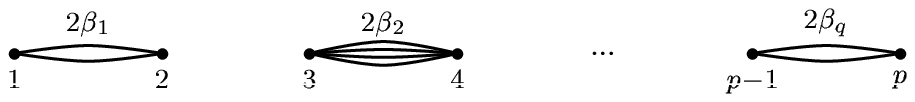}\end{gathered}
\]
where $2\beta_k = 2\alpha_{2k-1} = 2\alpha_{2k}$ is the multiplicity
of the edge which connects the vertices $2k-1$ and $2k$. The
$k$-th component of the reduced diagram $\gamma^*$ contributes by
the factor $|\rho (t_{2k-1}, t_{2k})|^{2\beta_k}$, so that the
value of diagram $\gamma$ is $\prod_{k=1}^q |\rho (t_{2k-1},
t_{2k})|^{2\beta_k}$. Integrating this over $T^p$ against the
measure $\prod_{l=1}^p \Theta (t_l) d\mu (t_l)$, we get
\begin{equation}\label{eq1.11}
\prod_{k=1}^q \iint_{T^2} |\rho (s,t)|^{2\beta_k} \Theta (s)
\Theta (t) \, d\mu (s) d\mu (t)\,.
\end{equation}

Now we need to do some counting. Each reduced diagram $\gamma^*$
is defined by the choice of the partition, like \eqref{eq1.10}, and the
multiplicities $\beta_1, \dots , \beta_q$ of the edges. Since the
glueing procedure is not one-to-one, each reduced diagram has a
``multiplicity'' (that is, the number of regular diagrams which
give the same reduced diagram $\gamma^*$ after glueing) depending on
$\beta_1, \dots , \beta_q$ which also should be taken into account.

Notice that if we started with a different partition of the
set $\{1, 2, \dots , p\}$ into disjoint pairs, then anyway we would
finish with the same answer \eqref{eq1.11}. Therefore, we must multiply
the expression \eqref{eq1.11} by the number of partitions of the set $\{1,
2, \dots , p\}$ into disjoint pairs. The number of such partitions
equals $\Ex \xi^p$ where $\xi \sim \cN(0,1)$ \cite[Remark~1.29]{Ja}.

Next, let us count the ``multiplicity'' of the reduced diagram. If
we fix a regular diagram $\gamma $, then all other regular
diagrams having the same reduced diagram $\gamma^*$ can be
obtained from $\gamma $ by $p$ independent permutations: we can
permute $\alpha_1$ vertices labeled by $1$, then $\alpha_2$
vertices labeled by $2$ and so on. Since these permutations are
independent, the total multiplicity of the reduced diagram equals
\begin{equation}\label{eq1.12}
\alpha_1! \alpha_2! \, \dots \, \alpha_p! = (\beta_1!)^2
(\beta_2!)^2\, \dots \, (\beta_q!)^2,
\end{equation}
and we need to put this factor before the product \eqref{eq1.11} when we
summate over reduced diagrams; i.e. over all possible
choices of the numbers $\beta_1, \dots , \beta_q \ge \alpha_0$.

Combining all pieces, we get
\begin{multline*}
\langle \text{regular} \rangle = \sum_{\beta_1,\dots ,\beta_q \ge
 \alpha_0}
 \frac{c_{2\beta_1}^2\dots c_{2\beta_q}^2}{(\beta_1!)^2\dots (\beta_q!)^2}\,
 (\beta_1 !)^2\, \dots \, (\beta_q!)^2 \, \Ex \xi^p \langle
 \text{product \eqref{eq1.11}} \rangle \\
= \left( \Ex \xi^p\right) \left( \sum_{\beta \ge \alpha_0}
 c_{2\beta}^2 \iint_{T^2} |\rho(s,t)|^{2\beta} \Theta (s)
 \Theta (t) \, d\mu (s) d\mu (t) \right)^q
\stackrel{\eqref{eq1.7}}= \left( \Ex \xi^p \right)\left( \Ex
Z_n^2 \right)^{p/2} \, .
\end{multline*}
Therefore,
\[
\Ex Z_n^p = \left( \Ex \xi^p \right)\left( \Ex
Z_n^2 \right)^{p/2} + \langle \mbox{irregular} \rangle\,,
\]
and we need to show that the contribution of irregular diagrams is
negligible with respect to the main term $\left( \Ex Z_n^2
\right)^{p/2}$ when $n\to\infty$.

Since $\phi $ is a polynomial, there are only finitely many
irregular diagrams which enter expression \eqref{eq1.8} for $\Ex
Z_n^p$ after plugging in the diagram formula \eqref{eq1.9}. We shall show
that if the diagram $\gamma $ is irregular, then
\begin{equation}\label{eq1.13}
\idotsint_{T^p} |V_\gamma (t_1, t_2, \dots , t_p)|\, d\mu
(t_1) \dots  d\mu (t_p) = o\left(\left( \sup_{s\in T} \int_T
|\rho_n(s,t)|\, d\mu (t)  \right)^{p/2}\right),
\end{equation}
for $n\to\infty$. Then due to assumption \eqref{eq1.3} (and expression
\eqref{eq1.7} for $\Ex Z_n^2$)
\[
\idotsint_{T^p} |V_\gamma (t_1, t_2, \dots , t_p)|\, d\mu
(t_1) \dots  d\mu (t_p) = o\left(\left( \Ex Z_n^2
\right)^{p/2}\right),
\]
and
\[
\langle \mbox{irregular} \rangle = o\left(\left( \Ex Z_n^2
\right)^{p/2}\right), \qquad n\to\infty\,,
\]
which finishes the proof of Theorem~1.2 in the case of general
radial functions  up to the proof of \eqref{eq1.13}.

\medskip\par\noindent{\em Proof of \eqref{eq1.13}:} First, we make another
reduction of the diagram and define a $p$-vertex graph with simple
edges which couple the vertices $i$ and $j$ if and only if at
least one of the pairs $(i, \overline{j})$ or $(\overline{j}, i)$
was coupled in the original diagrams (without taking into account
the multiplicities of the original coupling). We denote the
reduced diagram by $\gamma^{**}$. For example:
\[
\begin{gathered}\includegraphics{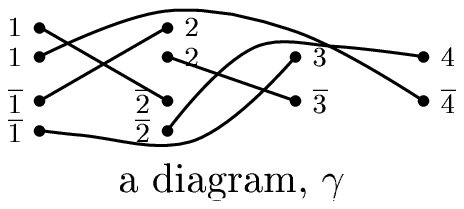}\end{gathered}\qquad\qquad
\begin{gathered}\includegraphics{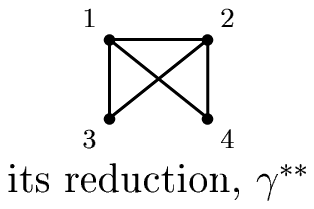}\end{gathered}
\]

Then
\begin{equation}\label{eq1.14}
| V_\gamma (t_1, \dots , t_p) | \le \prod_{(i,j)\in e(\gamma^{**})} |\rho
(t_i, t_j)|
\end{equation}
where the product is taken over all edges of $\gamma^{**}$. We
have to estimate from above the integral of $|V_\gamma|$ over
$T^p$. Replacing $|V_\gamma|$ by its upper bound \eqref{eq1.14}, we obtain
the integral which factorizes into the product of integrals
described by connected components of the diagram $\gamma^{**}$.

Let us start with one $m$-vertex component of the diagram
$\gamma^{**}$. The component can be a complicated graph --- anyway,
we can always turn this graph into a tree with $m$ vertices by
deleting some edges (this procedure only increases the integral we
are estimating). Having a tree, we choose a vertex belonging to
only one edge and integrate it out, which gives the factor $\left(
\sup_{s\in T} \int_T |\rho_n (s,t)|\, d\mu (t) \right)$ and the
rest of the tree which is a new tree with $m-1$ vertices. By
induction, any $m$-vertex tree describes the integral which does
not exceed
\[
\left( \sup_{s\in T} \int_T |\rho_n (s,t)|\, d\mu (t) \right)^{m-1}.
\]

Now, suppose the reduced diagram $\gamma^{**}$ has $k$ connected
components and the $i$-th component has $m_i$
vertices.\footnote{%
 Observe that $ m_1 + \dots + m_k = p $.}
Then
the right-hand side of \eqref{eq1.14} integrated over $T^p$ does not
exceed
\[
\left( \sup_{s\in T} \int_T |\rho_n (s,t)|\, d\mu (t)
\right)^{(m_1-1)+\dots+(m_k-1)} = \left( \sup_{s\in T} \int_T
|\rho_n (s,t)|\, d\mu (t) \right)^{p-k}.
\]
Since the diagram $\gamma $ is irregular, $k<p/2$ and we get
\eqref{eq1.13}.

\subsection{The last step}

It remains to explain why in the case when the
function $\phi$ is increasing, it suffices to assume that condition
\eqref{eq1.3} holds only for $\alpha = 1$.

In the proof given above, condition \eqref{eq1.3} was used only in the estimate
\[
\sup_{s\in T} \int_T |\rho_n (s,t)|\, d\mu (t) \le
\frac{\Ex Z_n^2}{c_{2\alpha_0}^2 \eta (\alpha_0)}
\]
where $\alpha_0\ge 1$ is the minimal positive integer such that
$c_{2\alpha_0}\ne 0$, and $\eta (\alpha)$ is the left-hand side of \eqref{eq1.3};
i.e., {\em we need condition \eqref{eq1.3} only with $\alpha = \alpha_0$}. If the
function $\phi$ increases, then $\alpha_0=1$. Indeed, it is easy to
find that $ \Wi{ |\zeta|^2 } = |\zeta|^2 -1 $. Then
\[
c_2 = \int_0^\infty \phi (r) (r^2-1) r e^{-r^2/2}\, dr
= \int_0^\infty \left( \phi (r) - \phi (1) \right) \left( r^2-1
\right) re^{-r^2/2}\, dr >0\,,
\]
concluding the proof of Theorem 1.2.


%% file: sect3.tex
\section{Asymptotic normality for chaotic analytic zero points}

\subsection{From chaotic analytic zeroes to non-linear
 functionals of Gaussian processes}

Recall the random objects defined in the
Introduction: the Gaussian analytic functions $\psi_L (z)$ (see
\eqref{elliptic}, \eqref{flat} and \eqref{hyperbolic}), the random measures $n_{\psi_L}$ counting
their zeroes, and the linear statistics $ Z_L(h) $. Since
\[
2\pi dn_\psi = \Delta \log|\psi| \, dm = \Delta^* \log|\psi| \, dm^*
\]
where $ dm^* $ and $ \Delta^* $ are the invariant measure and the
invariant Laplacian on $ \M $, we have
\[
Z_L(h) = \frac1{2\pi} \int_{\M} \log|\psi_L(z)| (\Delta^*
h)(z) \, dm^*(z)
\]
obtaining in each of the three cases a family of non-linear
functionals of the complex-valued Gaussian process $ \psi_L $
defined on $ \M $. We normalize the process $ \psi $
putting
\begin{equation}\label{eq3.1}
w_L(z) = \frac{\psi_L (z)}{\sqrt{\mathbb E|\psi_L(z)|^2}} =
\begin{cases}
 \psi_L (z) (1+|z|^2)^{-L/2} & \text{(elliptic)}, \\
 \psi_L(z) e^{-L|z|^2/2} & \text{(flat)}, \\
\psi_L(z) (1-|z|^2)^{L/2} & \text{(hyperbolic)}.
\end{cases}
\end{equation}
In the flat and hyperbolic cases we use a large real parameter
instead of the large integer parameter used in Section~2.
Note that the
(non-Gaussian) random process $\log |w(z)|$ is stationary (that
is, its distribution is invariant with respect to the
corresponding group of isometries), though the Gaussian process
$w(z)$ is not (because of the phase multipliers in \eqref{eq0.35}). In all
three cases, $w(z)\sim \mathcal N_{\mathbb C}(0,1)$ for each $z$.
We obtain
\[
Z_L(h)-\mathbb E Z_L(h) = \frac1{2\pi} \int_{\mathbb M} \log|w_L|
\Delta^* h\, dm^*\,.
\]
Indeed, the left-hand side has zero expectations since $\Ex
\log |w_L(z)| = \const $, and the invariant Laplacians of the
test functions are orthogonal to the constants. It remains to
check conditions \eqref{eq1.3} (with $\alpha=1$) and \eqref{eq1.4} of
Theorem~\ref{2.2},
and to compute the asymptotics of the variance of the random
variables $ Z_L(h) $.

\subsection{Checking conditions \eqref{eq1.3} and \eqref{eq1.4}}

Let
\[
\rho_L (z_1, z_2) = \mathbb E \left( w_L(z_1)
\overline{w_L(z_2)}\right) = \frac{\mathbb E \left( \psi_L(z_1)
\overline{\psi_L(z_2)} \right)}{\sqrt{\mathbb E |\psi_L(z_1)|^2
\cdot \mathbb E |\psi_L(z_2)|^2 }}\,.
\]
The function of two variables, $ |\rho_L (z_1,z_2)| $, reduces to
a function of one variable, $ |\rho_L (z,0)| $, by $ G
$-invariance, $ |\rho_L (z_1,z_2)| = |\rho_L \( g(z_1), g(z_2) \)|
$. Explicit formulas for $ \rho_L (z,0) = |\rho_L(z,0)| $ follow
from \eqref{elliptic}, \eqref{flat}, \eqref{hyperbolic} and
\eqref{eq3.1}; namely,
\begin{align*}
\rho_L(z,0) &= (1+|z|^2)^{-L/2} &&\text{(elliptic)}, \\
\rho_L(z,0) &= \exp \( - \frac12 L |z|^2 \) &&\text{(flat)}, \\
\rho_L(z,0) &= (1-|z|^2)^{L/2} &&\text{(hyperbolic)}.
\end{align*}

A straightforward inspection shows that in all three cases the
functions $\frac{L}{2\pi} \rho_L(z, 0)$ converge weakly to the
point mass at the origin:
\[
\lim_{L\to\infty} \frac{L}{2\pi} \int_\M \rho_L(z,0) \Theta(z)\,
dm^*(z) = \Theta(0)
\]
for any continuous test-function $\Theta$ (as usual, with a
compact support if $\M$ is non-compact).

By $G$-invariance, for any $\beta>0$, any continuous test-function
$ \Theta $ and any $z_2\in\M$,
\[
\lim_{L\to\infty} \frac{L\beta}{2\pi} \int_\M |\rho_L(z_1,
z_2)|^\beta \Theta(z_1) \, dm^*(z_1) = \Theta(z_2)
\]
(recall that $|\rho_L|^\beta = |\rho_{\beta L}|$ ). At last,
multiplying the both sides of this equation by $\Theta(z_2)$ and
integrating by $z_2$, we get
\begin{equation}\label{eq3.a}
\lim_{L\to\infty} \frac{L\beta}{2\pi} \iint_{\M^2} |\rho_L(z_1,
z_2)|^\beta \Theta(z_1) \Theta(z_2)\, dm^*(z_1) dm^*(z_2) =
||\Theta||^2_{L^2(m^*)}\,.
\end{equation}

Now, conditions \eqref{eq1.3} and \eqref{eq1.4} become evident:
for $L\to\infty$,
\[
\sup_{z_2\in \M} \int_\M |\rho_L(z_1, z_2)|\, dm^*(z_1) = \int_\M
|\rho_L(z, 0|\, dm^*(z) \sim \frac{2\pi}{L}\,,
\]
this gives us \eqref{eq1.4}. To get \eqref{eq1.3}, observe that
the double integral in the numerator of \eqref{eq1.3} has the same
decay
\[
\iint_{\M^2} |\rho_L(z_1, z_2)|^2 \Theta(z_1) \Theta (z_2) \,
dm^*(z_1) dm^*(z_2) \sim \frac{\pi}{L} \cdot
||\Theta||^2_{L^2(m^*)}\,.
\]

\subsection{Asymptotics of the variance}

We shall use formula \eqref{eq1.7} for the variance of $ Z(h) $,
\[
\Ex Z_L^2(h) = \frac1{4\pi^2} \sum_{\alpha\ge 1} c_{2\alpha}^2
\iint_{\M^2} | \rho_L (z_1,z_2) |^{2\alpha} \Delta^*h(z_1)
\Delta^*h(z_2) \, dm^*(z_1) \, dm^*(z_2) \, ,
\]
where $ c_{2\alpha} $ are defined by the expansion
\[
\log |\zeta| = \sum_{\alpha =0}^\infty \frac{ c_{2\alpha} }{
\alpha ! } :|\zeta|^{2\alpha}: \qquad \mbox{for } \zeta \sim
\mathcal N_\mathbb C(0,1)
\]
(cf.~\eqref{eq1.5}). Denoting
\[
\kappa = \frac1{4\pi} \sum_{\alpha \ge 1}
\frac{c_{2\alpha}^2}\alpha\,,
\]
we reduce \eqref{0.8L} to relation \eqref{eq3.a} with
$\beta=2\alpha$:
\[
\Ex Z_L^2(h) = \frac{1+o(1)}{4\pi^2} \sum_{\alpha\ge 1}
c_{2\alpha}^2 \cdot \frac{\pi}{L\alpha} \cdot
||\Delta^*h||^2_{L^2(m^*)} = \frac{\kappa + o(1)}{L} \cdot
||\Delta^*h||^2_{L^2(m^*)}\,.
\]